\documentclass[11pt]{article}
\usepackage{amsmath}
\usepackage{latexsym}
\usepackage{amsfonts}
\usepackage{amssymb}
\usepackage{theorem}

\evensidemargin  \oddsidemargin

\numberwithin{equation}{section}
\newtheorem{theorem}{Theorem}[section]

\newtheorem{lemma}[theorem]{Lemma}

{\theorembodyfont{\rmfamily} }

\newcommand{\bke}[1]{\left ( #1 \right )}
\newcommand{\bkt}[1]{\left [ #1 \right ]}
\newcommand{\bket}[1]{\left \{ #1 \right \}}
\newcommand{\norm}[1]{\left \| #1 \right \|}
\newcommand{\bka}[1]{{\langle #1 \rangle}}
\newcommand{\abs}[1]{\left | #1 \right |}
\newcommand{\e}{\varepsilon}

\newcommand{\la}{\lambda}
\newcommand{\al}{\alpha}
\newcommand{\pd}{\partial}
\newcommand{\ph}{\phi}
\newcommand{\de}{\delta}
\newcommand{\si}{\sigma}
\newcommand{\ga}{\gamma}
\renewcommand{\th}{\theta}
\newcommand{\R}{\mathbb{R}}

\newcommand{\Z}{\mathbb{Z}}

\newcommand{\lec}{{\ \lesssim \ }}

\newcommand{\sgn}{\mathop{\mathrm{sgn}}}
\newcommand{\les}{\lesssim}
\newcommand{\spt}{\mathop{\mathrm{supp}}}
\newcommand{\I}{{\cal I}}
\newcommand{\G}{{\cal G}}
\newcommand{\sumk}{{\textstyle \sum_{k}}}
\newcommand{\qed}{\hfill$\Box$}
\newcommand{\cH}{{\cal H}}

\begin{document}
\title{Stability in $H^{1/2}$ of the sum of $K$ solitons for the
Benjamin-Ono equation}

\author{ Stephen Gustafson,\quad Hideo Takaoka, \quad Tai-Peng Tsai}

\date{} 

\maketitle

\begin{abstract}
This note proves the orbital stability in the energy space $H^{1/2}$
of the sum of widely-spaced $1$-solitons for the Benjamin-Ono
equation, with speeds arranged so as to avoid collisions.
\end{abstract}

\section{Introduction}
In this article we study the stability problem of the sum of $K$ solitons
for the Benjamin-Ono (BO) equation for 
$u(t,x) : \R^+ \times \R \to \R$:
\begin{equation}\label{BO1} 
  u_t = -({\cal H} \pd_x u + u^2)_x,
\end{equation}
where ${\cal H}$ is the Hilbert transform operator defined by
$$
{\cal H} f(x) =\mbox{p.v.}\frac 1 \pi\int_\R \frac {f(y)}{x-y}\,dy.
$$
Alternatively, if we denote $D =
\sqrt{-\pd_x^2}$, we have ${\cal H} \pd_x = -D$ and we can rewrite 
the Cauchy problem for~\eqref{BO1} as\footnote{The Fourier transform is given by
\[
\hat f (\xi) = \frac 1{\sqrt{2\pi}} \int_\R f(x) e^{-i x\xi}\, dx,\quad
f (x) = \frac 1{\sqrt{2\pi}} \int_\R \hat f(\xi) e^{i x\xi}\, d\xi
\]
so that
$\widehat{\pd_xf}(\xi) = -i\xi \hat f (\xi)$, $\widehat{{\cal H}f}(\xi)
= -i\sgn(\xi) \hat f (\xi)$ 
and $\widehat{Df}(\xi) = |\xi| \hat f (\xi)$.}
\begin{equation}\label{BO2}
\begin{split}
  &u_t = ( Du - u^2)_x \\
  &u(0,x) = u_0(x).
\end{split}
\end{equation}
This equation is a model for one-dimensional long waves in deep stratified fluids (\cite{Ben,Ono}).

The Benjamin-Ono equation is completely integrable and has infinitely
many conserved quantities (\cite{KLM1,KLM2}).  Two of them are the
$L^2$ mass
$$
N(u) = \frac 12 \int_\R u^2 \,dx,
$$
and the energy
$$
E(u) = \int_\R  \frac 12 u D u - \frac 13 u^3 \, dx.
$$ 
The energy space, where $E(u)$ is defined, is $H^{1/2}(\R)$.  The
existence of global {\it weak} solutions $u \in
C([0,\infty);H^{1/2}(\R)) \cap C^1((0,\infty);H^{-3/2}(\R))$ to
\eqref{BO2} with energy space initial data $u(0,x) = u_0(x) \in
H^{1/2}(\R)$ was shown by J. C. Saut \cite{Saut} (see also the paper
of J. Ginibre and G. Velo \cite{GV}).  For the {\it strong}
$H^s$-solution, A. Ionescu and C. E. Kenig \cite{IK} established
global well-posedness for $s\ge 0$ (see also the paper of T. Tao
\cite{Tao}).  This solution conserves the functional $N(u)$ (and
$E(u)$ when $s\ge 1/2$).

The Benjamin-Ono equation admits ``$K$-soliton'' solutions \cite{Joseph}.
The 1-solitons are of the form
$$
u(t,x) = Q_c(x-ct-x_0), \quad (c>0, \ x_0 \in \R)
$$
where
\begin{equation}
\label{Q.formula}
Q_c(x) = c Q(cx), \quad Q(x) = \frac 2 {1+x^2}.  
\end{equation}
They satisfy
\begin{equation}
{\cal H}\pd_x Q_c + Q_c^2 = c Q_c,
\end{equation}
which can be verified by using
$\hat Q(\xi) = \sqrt{2 \pi }\, e^{-|\xi|}$.
By the explicit form~\eqref{Q.formula}, we have
\begin{equation}\label{eq:QDQ}
\int Q^2 = 2\pi, \quad \int Q^3 = 3 \pi, \quad (Q,DQ)= \int Q(Q^2-Q) =
\pi.
\end{equation}
By rescaling,
\begin{equation*}
N(Q_c) = c N(Q) = \pi c, \quad
E(Q_c) = c^2 E(Q) = -\frac \pi 2 c^2.
\end{equation*}

The orbital (i.e. up to translations) stability of the $1$-soliton
in the energy norm ($H^{1/2}$) was established in \cite{BBSSB}.
See \cite{Be,Bona} for earlier stability results.
Here we address the stability of the sum of widely-spaced 
$1$-solitons, with speeds arranged so as to avoid collisions.
Our main result is the following theorem.
\begin{theorem} [Orbital stability of the sum of $K$ solitons]
\label{TH1-1}
Let $0< c_1^0 < \cdots < c_K^0$.
There exist $L_0, A_0, \al_0>0$ and $\th_0 \in (0,1)$ such that for any
$u_0 \in H^{1/2}(\R)$, $L>L_0$, and $0<\al<\al_0$, if
\begin{equation}
\label{data}
  \| u_0 - \sum_{k=1}^K Q_{c_k^0}(\cdot -x_k^0)\|_{H^{1/2}(\R)} \le \al
\end{equation}
for some  $x_k^0$ satisfying
\begin{equation}
\label{spacing}
  x_{k+1}^0 - x_k^0 >L, \quad (k=1, \ldots K-1),
\end{equation}
then there exist $C^1$-functions $x_k(t)$,
$k=1,\ldots, K$, such that the solution of~\eqref{BO2} satisfies
\[
\| u(t) - \sum_{k=1}^K Q_{c_k^0}(\cdot -x_k(t))\|_{H^{1/2}(\R)}
\le A_0 (\al + L^{-\th_0}), \quad \forall t >0.
\]
Moreover,
\[
|\dot x_k(t) - c_k^0 | \le  A_0 (\al + L^{-\th_0}), \quad \forall t >0.
\]
\end{theorem}

Integrable systems techniques (in particular higher conservation
laws) have been used to establish the stability of exact $K$-soliton
solutions (see \cite{MS} for KdV, and \cite{NL} for BO $2$-solitons)
against perturbations which are small in (necessarily) higher 
Sobolev norms. Here we are considering
a different problem: stability of sums of $1$-solitons (configurations
which are not themselves solutions) in the {\it energy space}.
Results of this type were obtained for KdV-type equations and 
NLS equations in \cite{MMT1,DM,DMo} and \cite{MMT2}, respectively.
Our approach follows that of \cite{MMT1} for gKdV, which adds to the
energy method of Weinstein \cite{W86} for the one soliton case, the
monotonicity property of the $L^2$-mass on the right of each soliton.
Here we encounter two new difficulties. Firstly, and most importantly,
the operator ${\cal H}$ is non-local, necessitating commutator estimates.
Secondly, the decay of the soliton $Q(x)$ is only algebraic,
meaning the error estimates are more delicate.
In particular, we use cut-off functions whose supports expand
sublinearly at the rate $O(t^\gamma)$, $2/3<\gamma<1$, similar
to \cite{MMT2}.

After the paper was completed, we learned that C.~E.~Kenig and Y.~Martel
\cite{KM} have obtained a similar result independently and
simultaneously.

\section{The stability proof}

Here we prove Theorem~\ref{TH1-1} using a series of Lemmas
whose proofs are given in section 3.

So we begin by fixing speeds $0 < c_1^0 < \cdots < c_K^0$,
and we suppose 
$u \in C([0,\infty);H^{1/2}(\R)) \cap C^1((0,\infty);H^{-3/2}(\R))$
solves~\eqref{BO2} with initial data satisfying~\eqref{data}
and~\eqref{spacing} for $\alpha < \alpha_0$ and $L > L_0$,
where $\alpha_0\ll 1$ and $L_0 \gg 1$ will be determined (depending only on
the speeds $\{ c_k^0 \}$) in the course of the proof.

\subsection{Decomposition of the solution}

Set
\begin{equation}
\label{close}
  T = T(\alpha,L) := \sup \; \Big\{ \; t > 0 \; | \; \sup_{0 \leq s \leq t}
  \; \inf_{y_j > y_{j-1} + L/2}
  \Big\| u(s,.)-\sum_{j=1}^K Q_{c^0_j}(.-y_j)\Big\|_{H^{1/2}}
  < \sqrt{\alpha} \; \Big\}.
\end{equation}
If we take $\alpha < 1$, then since $u \in C([0,\infty);H^{1/2})$,
we have $T > 0$. In what follows, we will estimate on the time
interval $[0,T]$, and in the end conclude (provided $\alpha$ 
sufficiently small, $L$ sufficiently large) that $T = \infty$.

The first step is a decomposition of the solution.
\begin{lemma}[Decomposition of the solution]\label{L2.1}
There exist $L_1 > 0$, $\alpha_1 > 0$, and $A_1>0$
such that if $\alpha < \alpha_1$and $L > L_1$, then
there exist unique $C^1$-functions
$c_j:[0,T] \to (0,+\infty)$, $x_j:[0,T] \to \R$, such that
\begin{equation}\label{defofeps}
  u(t,x) = \sum_{j=1}^K R_j(t,x) + \e(t,x)
  \quad \hbox{where} \quad R_{j}(t,x) := Q_{c_j(t)}(x-x_j(t)),
\end{equation}
where $\e(t,x)$ satisfies the orthogonality conditions
\begin{equation}\label{ortho}
\forall j, \forall t\in [0,T],\quad \int R_j(t,\cdot) \e(t,\cdot)=
\int (R_j(t,\cdot))_x \e(t,\cdot)=0.
\end{equation}
Moreover, 
\begin{equation}
\label{initial}
  \| \e(0,\cdot) \|_{H^{1/2}} + \sumk |x_k(0) - x_k^0| 
  + \sumk |c_k(0) - c_k^0| \leq A_1 \alpha, 
\end{equation}
and for all $t \in [0,T]$, 
\begin{equation}\label{largeness}
  x_k(t) - x_{k-1}(t) \geq L/2,
  \quad\quad k=2,\ldots,K,
\end{equation}
\begin{equation}\label{smallness1}
  \| \e(t,\cdot) \|_{H^{1/2}} + \sum_{j=1}^K |c_j(t) - c_j^0|
  \le A_1 \sqrt{\alpha},
\end{equation}
\begin{equation}\label{smallness2}
  \sum_{j=1}^K |\dot x_j(t)-c^0_j| + |\dot{c_j}(t)|
  \le A_1 (\sqrt{\alpha} + L^{-2}).
\end{equation}
\end{lemma}

We will use $\norm{\e(t,\cdot)}_{H^{1/2}} \le 1$ in the rest of the proof.

\subsection{Almost monotonicity of local mass}

The size of the remainder $\e(t,x)$ will be controlled by
an ``almost monotone'' Lyapunov functional which we now construct.
Fix
\[
  \gamma \in (2/3,1),
\]
and a nonnegative
$\zeta(x)\in C^2(\R)$ so that $\zeta(x)=1$ for $x >1 $,
$\zeta(x)=0$ for $x<0$, and $\sqrt{\zeta_x} \in C^1$. Set 
\[
  \bar x_k^0 := \frac{x_{k-1}(0)+x_k(0)}{2}, \quad\quad
  \sigma_k := \frac{c_{k-1}^0 + c_k^0}{2}, \quad\quad
  k = 2,\cdots,K,
\]
$\psi_1(t,x) \equiv 1$, and for $k=2,\ldots,K$,
\begin{equation}
\label{psidef}
  \psi_k(t,x) := \zeta(y_k), \quad\quad
  y_k := \frac{x - \bar x_k^0 - \si_k t}{(b+t)^{\gamma}},
\end{equation}
with
$b := \left( \frac{L}{16} \right)^{1/\gamma}$,
and, finally, set for $k=1,\ldots,K$,
\[
  \I_k(t) := \frac 12 \int_\R \psi_k(t,x) u(t,x)^2\, dx,
\]
which, roughly speaking, measures the $L^2$ mass to the right of the
$k$-th soliton.

Setting
\[
  d_k := c_k(0) - c_{k-1}(0), \quad (k=2,\dots, K); 
  \quad d_1 = c_1(0),
\]
the Lyapunov function we will use is
\begin{equation}
\label{Gdef}
  \G(t) := E(u(t)) + \sum_{k=1}^K d_k \I_k(t) .
\end{equation}
Note that $E(u(t))=E(u_0)$ by energy conservation.
The ``almost monotonicity'' of this functional comes from
the following key estimate.
\begin{lemma}[Almost monotonicity of mass on the right of each
    soliton]\label{L2.4}
Under the decomposition in \eqref{defofeps}, there is $C_2 > 0$ such
that
\[
\I_k(t) - \I_k(0) \le C_2 L^{\frac{1}{\gamma}-\frac{3}{2}} + C_2
  L^{1-\frac{1}{\gamma}} \sup_{t' \in [0,t]} \norm{\e(t')}_{L^2}^2 .
\]
\end{lemma}

In light of~\eqref{Gdef}, this lemma implies the estimate
\begin{equation}
\label{en}
  \G(t) \le \G(0) + CL^{\frac{1}{\gamma}-\frac{3}{2}}  
  +C L^{1-\frac{1}{\gamma}} 
  \sup_{t'\in [0,t]} \norm{\e(t')}_{L^2}^2.
\end{equation}

\subsection{Decomposition of the energy}

As above, set
$R_k := Q_{c_k(t)}(x-x_k(t))$, $R := \sum_{k=1}^K R_k$,
and define
\[
  \phi_k(t,x) := \psi_k(t,x) - \psi_{k+1}(t,x),
  \quad k=1,\ldots,K-1, \quad \phi_K(t,x) = \psi_K(t,x)
\]
(so $\phi_k$ is localized near the $k$-th soliton), and
the (time-dependent) operator
\[
  H_K := D - 2R + \sum_{k=1}^K c_k(t) \phi_k.
\]   
The functional $\G$ can be expanded as follows.
\begin{lemma}
[Energy decomposition]
\label{energy}
There is $C_3 > 0$ such that
\[
\begin{split}
  \Big| \G(t) -  \big\{
  \sum_k [ E(R_k) + & c_k(0) N(R_k) ]
  + \frac{1}{2} (\e(t), H_K \e(t)) \big\} \Big| \\
  &\leq C_3 \Big( L^{-2} + \|\e(t)\|_{L^2} 
  \sum_k |c_k(0) - c_k(t)| + \|\e(t)\|_{H^{1/2}}^3 \Big).
\end{split}
\]                               
\end{lemma}

We also need:
\begin{lemma}\label{L.Fdiff}
Let $F_c(u) = E(u) + c N(u)$. We have for some $C>0$ and
$c$ close to $c^0$ that
\[
  0 \le F_{c^0}(Q_{c^0}) - F_{c^0}(Q_c) \le C (c-c^0)^2.
\]
\end{lemma}

Combining equation \eqref{en}, Lemmas \ref{energy} and \ref{L.Fdiff} yields
\begin{equation}
\label{en3}
\begin{split}
  (\e(t),H_K\e(t)) &\le C \Big[ \sumk |c_k(0) - c_k(t)|\norm{\e(t)}_{L^2}
  + \norm{\e(t)}_{H^{1/2}}^3 + \sumk |c_k(0) - c_k(t)|^2 \\
  &\quad + \norm{\e(0)}_{H^{1/2}}^2 
  + L^{\frac{1}{\gamma}-\frac{3}{2}}
  + L^{1-\frac{1}{\gamma}} 
  \sup_{t'\in [0,t]}\norm{\e(t')}_{L^2}^2 \Big].
\end{split}
\end{equation}

Next we need quadratic control of $c_k(t)-c_k(0)$.
\begin{lemma}[Quadratic control of speed change]
\label{lem:c}
\[
\sumk |c_k(t) - c_k(0)| \lec L^{\frac{1}{\gamma}-\frac{3}{2}} +
  L^{1-\frac{1}{\gamma}} \sup_{0 \leq \tau \leq t}
  \norm{\e(\tau)}_{L^2}^2 + \norm{\e(t)}_{H^{1/2}}^2 +
  \norm{\e(0)}_{H^{1/2}}^2.
\]
\end{lemma}

Combining this lemma with~\eqref{en3} and setting 
$\th_0 = \frac 12(\frac32-\frac 1 \gamma)$ yields
\begin{equation}
\label{en4}
  (\e(t),H_K\e(t)) \lec 
  \norm{\e(t)}_{H^{1/2}}^3 
  + \norm{\e(0)}_{H^{1/2}}^2 +
  L^{-2\th_0} + L^{1-\frac{1}{\gamma}} 
  \sup_{t'\in [0,t]}\norm{\e(t')}_{L^2}^2.
\end{equation}

\subsection{Lower bound on quadratic form
and completion of the proof}

We want to use the quadratic form $(\e,H_K \e)$ to control
$\norm{\e}_{H^{1/2}}^2$, as is done for one-soliton
stability.
Here we need a $K$-soliton version of this.
\begin{lemma}[Positivity of the quadratic form] \label{TH3-2}
There exist $L_2, \ga_2 > 0$ such that if $L > L_2$,  then
\[
  \ga_2 \norm{\e }_{H^{1/2}}^2 \le (\e , H_K \e ).
\]
\end{lemma}

Combining this lemma with~\eqref{en4} gives
\[
  \norm{\e(t)}_{H^{1/2}}^2 \le C \Big[  \| \e(t) \|_{H^{1/2}}^3 +
  \norm{\e(0)}_{H^{1/2}}^2 + L^{-2\th_0} +
  L^{1-\frac{1}{\gamma}} \sup_{t'\in [0,t]}\norm{\e(t')}_{L^2}^2
  \Big].  
\]
So using~\eqref{initial}, this estimate implies, 
for $\alpha$ and $1/L$ sufficiently small,
that there is $A_0 > 0$ such that
\begin{equation}
\label{finalest}
  \sup_{t'\in[0,T]}\norm{\e(t')}_{H^{1/2}} 
  \le A_0 (\al + L^{-\theta_0}).
\end{equation}
Hence for $\alpha$ and $1/L$ sufficiently small, 
we conclude $T = \infty$, $x_k(t)$ and $c_k(t)$ exist
for all time, and~\eqref{finalest} gives the main
estimate of the theorem. 
Finally, the last estimate of the theorem follows 
from~\eqref{diff} and~\eqref{diff2} in the proof
of Lemma~\ref{L2.1}.   \hfill $\square$

\section{Proofs of lemmas}
In this section, we shall prove lemmas mentioned in section 2.
\subsection{Decomposition of the solution}
\noindent
{\it Proof of Lemma~\ref{L2.1}.}
The existence of the functions $c_j(t)$ and $x_j(t)$ is 
established through the implicit function theorem
applied to the map
$F : H^{-3/2}(\R) \times \R^K \times (\R^+)^{K} \to \R^{2K}$
defined by
\[
  F(u,{\bf y},{\bf c}) := 
  \big( ( {\bf R}, u - R ), ( {\bf R}_x, u - R ) \big)
\]
where $R(x) = \sum_{j=1}^K R_j(x)$ with $R_j(x) = Q_{{c_j}}(x-y_j)$,
and boldface denotes $K$-vectors, e.g.,~${\bf y} = (y_1,\ldots,y_K)$
and ${\bf R}(x) := (R_1(x),\ldots,R_K(x))$.  Here the inner product
indicates $H^{3/2} - H^{-3/2}$ pairing.  $F$ is easily seen to be a
$C^1$ map (note it is affine in $u$).  For any ${\bf y}$ and (bounded)
${\bf c}$, $F(R,{\bf y},{\bf c}) = ({\bf 0}, {\bf 0})$, and as a $2K
\times 2K$ matrix,
\begin{equation}
\label{invert}
  D_{{\bf y}, {\bf c}} F(R,{\bf y}, {\bf c}) = 
  \pi \left( \begin{array}{cc} 
  0 & -Id \\
   diag(c_j^3) & 0 \end{array} \right)
  + O \big((\min_{j \not= k} |y_j - y_k|)^{-2}\big)
\end{equation}
is invertible, provided 
$\min_{j \not= k} |y_j - y_k| >  L_1/2$ ($L_1$ a constant).
Thus there is $\alpha_1 > 0$ such that for any ${\bf y}$ 
satisfying this condition, for $u$ in an $H^{-3/2}$-ball about 
$\sum_{j=1}^K Q_{c_j^0}(x-y_j)$ of size 
$\beta \in (0,\sqrt{\alpha_1})$, 
there are unique $C^1(H^{-3/2};\R^K)$ functions 
${\bf x}(u)$ and ${\bf c}(u)$ so that
$F(u,{\bf x}(u),{\bf c}(u)) = 0$, with
\begin{equation}
\label{diff}
  |{\bf c}(u) - {\bf c}^0| + |{\bf x}(u) - {\bf y}| \lec 
  \beta .
\end{equation}
So using the condition~\eqref{close},
for $0 \leq t \leq T$, we take $\beta = \sqrt{\alpha}$
and set 
${\bf c}(t) := {\bf c}(u(t))$ and ${\bf x}(t) := {\bf x}(u(t))$.
Since $u \in C^1((0,\infty);H^{-3/2})$, 
$x_j(t)$ and $c_j(t)$ are $C^1$ functions of $t>0$. 
The equation $F(u,{\bf x}(t),{\bf c}(t)) \equiv 0$
is equivalent to the orthogonality conditions~\eqref{ortho}.
The estimates~\eqref{smallness1} 
follow from~\eqref{close} and~\eqref{diff}.
An equation for $\e(t,x)$ can be derived using \eqref{BO2} and
$(DR_k - R_k^2)_x - \pd_t R_k = (\dot x_k - c_k)\pd_x R_k - \dot c_k \pd_c R_k$:
\begin{equation}
  \pd_t \e = \pd_x( D\e - 2 R\e - \e^2 - \sum_{j \not = k}R_jR_k)
  +\sum_k (\dot x_k - c_k)\pd_x R_k - \dot c_k \pd_c R_k.
\end{equation}
Computing $\frac d{dt}(R_k,\e)$ and $\frac d{dt}(\pd_x R_k,\e)$
in turn, and using~\eqref{invert} and~\eqref{smallness1}
yields
\begin{equation}
\label{diff2}
  |\dot {\bf c}| + |\dot {\bf x} - {\bf c}| \lec 
  \| \e \|_{H^{1/2}} + L^{-2} 
  \lec \sqrt{\alpha} + L^{-2}.
\end{equation}
This implies that $c(t)$ and $x(t)$ are $C^1$ up to $t=0$ and, together
with~\eqref{diff}, it gives \eqref{smallness2}.  Now $\alpha$ can be
taken sufficiently small, and $L$ sufficiently large, so
that~\eqref{data}-~\eqref{spacing}, together with~\eqref{diff} with
$\beta=\alpha$, imply~\eqref{initial}, which in turn implies that
$x_k(0) - x_{k-1}(0) \geq L/2$.  Finally~\eqref{largeness} follows
from this and~\eqref{smallness2} via
\[
  \frac{d}{dt}(x_k(t) - x_{k-1}(t)) 
  \geq c_k^0 - c_{k-1}^0 - A_1(\sqrt{\alpha} + L^{-2}) > 0
\]
for $\alpha$ sufficiently small, $L$ sufficiently large.
\qed
\subsection{Commutator estimates}
We have to deduce several estimates for commutators.  For two
operators $A$ and $B$, denote by $[A,B]=AB-BA$ their commutator.
\begin{lemma} \label{L2.2}
(i) Suppose $\chi \in C^1_c(\R)$. We have
\begin{equation} \label{L2.2-a}
\norm{[D^{1/2}, \chi] u }_{L^2(\R)} \les \norm{{|\xi|^{1/2}\hat \chi (\xi)}
}_{L^1(d\xi)} \cdot \norm{u}_{L^{2}}.
\end{equation}
(ii) Suppose $\phi \in B^{2-2\e}_{\infty,1}$ with
$0<\e<1/2$, then
\begin{equation} \label{L2.2-b}
\left|\int u_x[{\cal H},\phi]u_x\right|\lesssim
\norm{\phi}_{\dot B^{2-2\e}_{\infty,1}}
\|u\|_{H^{1/2}}^2.
\end{equation}
\end{lemma}

{\it Proof.}
(i) One can show $\abs{|p|^{1/2}- |p - \xi|^{1/2}}\lec
|\xi|^{1/2}$ 
by considering the two cases $|p| > 3 |\xi|$
and $|p| < 3 |\xi|$. Thus
\begin{align*}
& \abs{ D^{1/2}(u\chi) - (D^{1/2}u) \chi}_{L^2(dx)}
= \abs{\int [{|p|^{1/2}- |p - \xi|^{1/2}] \hat u(p - \xi)
\hat \chi (\xi)} d \xi}_{L^2(dp)}
\\
& \lec \abs{|\hat u| * |\xi|^{1/2} |\hat \chi| }_{L^2}
 \le |\hat u|_2 \cdot \abs{ |\xi|^{1/2} \hat \chi }_{L^1}.
\end{align*}

(ii) First assume $\phi \in C^2_c(\R)$. Let
$\Gamma=\{\xi_1+\xi_2+\xi_3=0\}$.  The integral is equal to
$$
\int u_x[{\cal H},\phi]u_x =i\int_\Gamma
\xi_1\xi_3 m(\xi) \hat{u}(\xi_1)\hat{\phi}(\xi_2)\hat{u}(\xi_3)
$$
where
$$
m(\xi)= \sgn(\xi_2+\xi_3) -\sgn(\xi_3).
$$
Decompose the integral into a  sum by Littlewood-Paley decomposition
\[
\sum_{N_1,N_2,N_3} \int_{\Gamma}\xi_1\xi_3 m(\xi)
\hat u_{N_1}(\xi_1)\hat{\phi}_{N_2}(\xi_2)\hat{u}_{N_3}(\xi_3)
\]
where $N_j$ are dyadic numbers, $N_j=2^k$ for $k \in \Z$.  

If $N_3 \gg N_2$, then $m(\xi)=0$. If $N_3 \les N_2$, then $N_1 \les
N_2$ on $\Gamma$.  Thus we may assume $N_1,N_3 \lec N_2$.  When
$\xi_1+\xi_2+\xi_3=0$, $m(\xi)=m_1(\xi)+m_3(\xi)$ where $m_j(\xi)=-
\sgn(\xi_j)$ is constant when $\xi_j \not =0$.  By multi-linear
estimates \cite[Theorem 1.1]{MTT}, we have
\[
\abs{\int_{\Gamma} m(\xi) 
 \xi_1 \hat u_{N_1}(\xi_1)\hat{\phi}_{N_2}(\xi_2) \xi_3\hat{u}_{N_3}(\xi_3)}
\le C \norm{\nabla u_{N_1}}_{2} 
\norm{\phi_{N_2}}_{\infty} \norm{\nabla u_{N_3}}_{2} .
\]
Thus
\begin{align*}
\abs{  \int u_x[{\cal H},\phi]u_x }
&\lec \sum_{N_1,N_3 \lec N_2} N_1 N_3 \norm{u_{N_1}}_{2} 
\norm{\phi_{N_2}}_{\infty} \norm{ u_{N_3}}_{2}
\\
&\lec \sum N_1^\e N_2^{2-2\e} N_3^\e \norm{u_{N_1}}_{2} 
\norm{\phi_{N_2}}_{\infty} \norm{ u_{N_3}}_{2} 
\\
&=  \norm{\phi}_{\dot B^{2-2\e}_{\infty,1}} \norm{u}_{\dot B^\e_{2,1}}^2 .
\end{align*}
Since $\norm{u}_{\dot B^\e_{2,1}} \lec \norm{u}_{ H^{1/2}}$
for $0<\e<\frac 12$, we have shown \eqref{L2.2-b} for $\phi \in C^2_c$.

For general $\phi \in B^{2-2\e}_{\infty,1}$, take $\eta_R(x) =
\eta(x/R)$ where $\eta(x)$ is a fixed smooth function which equals 1
for $|x|<1$ and 0 for $|x|>2$.  We have
\[
\norm{\phi \eta_R}_{\dot B^{2-2\e}_{\infty,1}} \lec \norm{\phi
}_{L^\infty} \norm{ \eta_R}_{\dot B^{2-2\e}_{\infty,1}}
+\norm{\eta_R }_{L^\infty} \norm{\phi }_{\dot B^{2-2\e}_{\infty,1}}.
\]
Sending $R$ to infinity in \eqref{L2.2-b} with the above estimate, we
get  \eqref{L2.2-b} for $\phi \in B^{2-2\e}_{\infty,1}$.
\qed

\begin{lemma} \label{L2.3}
Suppose $\chi(x) \in C^1_c(\R)$ and $\widehat{D^{1/2} \chi} \in L^1$.
For all $u \in H^{1/2}$,
\[
\int_\R u^3 \chi^2 \, dx \le C
\bke{\int_{\spt \chi} u^2 dx}^{1/2}
\int \bke{|D^{1/2}u|^2 \chi^2 + u^2\chi^2
+ u^2 \left\|\widehat{D^{1/2} \chi} \right\|_{L^1_\xi}^2 }\,dx.
\]
Here $C$ is a constant independent of $u$ and $\chi$.
\end{lemma}
\noindent
{\it Proof.}
First note the Gagliardo-Nirenberg inequality 
\begin{equation}
\int u^4\, dx \les \int |D^{1/2}u|^2\,dx \cdot
\int |u|^2\,dx  .
\end{equation}
This can be proved by first noting
\[
\|u\|_{4} \les \|\hat u\|_{4/3}  \le \|\bka{\xi}^{1/2} \hat u\|_2
\|\bka{\xi}^{-1/2}\|_{4} = C (\|D^{1/2}u\|_2 + \|u\|_2),
\]
and then rescaling with a minimizing scaling parameter. By H\"older
inequality and the above inequality,
\begin{align*}
\Big(\int u^3 \chi^2 dx \Big)^2   &\le \int_{\spt \chi} u^2 dx
\int (u \chi)^4 dx
\\
&\le \int_{\spt \chi} u^2 dx \int |D^{1/2}(u\chi)|^2\,dx 
\int u^2 \chi ^2 \,dx .
\end{align*}
By equation \eqref{L2.2-a}, we conclude
\[
\Big(\int u^3 \chi^2 dx \Big)^2   \lec \int_{\spt \chi} u^2 dx
\int u^2 \chi ^2 \,dx
\bke{ \int |D^{1/2}u|^2 \chi^2\,dx
+ \int u^2 dx \| |\xi|^{1/2} \hat \chi \|_{L^1}^2},
\]
from which the lemma follows.
\qed
\subsection{Almost monotonicity}
\noindent
{\it Proof of Lemma~\ref{L2.4}.} We may assume $u$ is smooth since the
general case follows from approximation.  We may assume $k \ge 2$
since $\I_1(t)$ is constant.  Denote $\psi = \psi_k(t,x) = \zeta(y_k)$
for simplicity of notation.  Note $\psi \in B^{2-}_{\infty,1}$ and
\begin{equation}
\label{eq2-4}
  \psi_x = (b+t)^{-\gamma} \zeta'(y_k) , \quad
  \spt \psi_x \subset \bar x_k^0 + \si_k t + [0, (b+t)^{\gamma}].
\end{equation}
Consider
\begin{align*}
\frac d{dt} \I_k(t)
&= \int - \psi u [\cH u_{x}+u^2]_x +\frac 12 u^2 \pd_t \psi \, dx
\\
& =  \int (\psi_x u + \psi u_x) \cH u_{x} + \frac 23 u^3 \psi_x
- \frac 12 u^2 \bke{ \si_k \psi_x + \zeta'(y_k)
 \frac \gamma {b+t} y_k} \, dx
\end{align*}

By $\cH\pd_x =- D$ and by Lemma \ref{L2.2} (i) with $\chi = \psi_x$, we
have
\begin{align*}
\int \psi_x u  \cH u_{x} &= - \int \psi_x |D^{1/2}u|^2 - \int (D^{1/2}u)
[D^{1/2},\psi_x] u
\\
& = - \int \psi_x |D^{1/2}u|^2
+ O(\norm{u}_{H^{1/2}}^2 \norm{ |\xi|^{1/2} \widehat{\psi_x}}_{L^1})
\end{align*}
Since $ \int \psi u_x  \cH u_{x} = -\int  u_x  \cH (\psi u_{x} )
= -  \int \psi u_x  \cH u_{x}  -\int  u_x  [\cH , \psi]u_x$,
by Lemma \ref{L2.2} (ii),
\[
\int \psi u_x  \cH u_{x} = - \frac 12 \int  u_x  [\cH , \psi]u_x
= O(\norm{u}_{H^{1/2}}^2  \|\psi\|_{\dot B^{2-2\e}_{\infty,1}}).
\]
Here we choose $\e \in (0,\frac 14)$. By Lemma \ref{L2.3} with $\psi_x
= \chi^2$,
\[
\int \frac 23 u^3 \psi_x \lec
\norm{u}_{L^{2}(\spt \psi_x)}\cdot
 \int (|D^{1/2}u|^2+u^2)\psi_x
+ u^2 \norm{ |\xi|^{1/2} {\cal F}( \sqrt{\psi_x})}_{L^1_\xi} ^2.
\]
Now by~\eqref{initial},~\eqref{smallness1}, and the definition of 
$\sigma_k$, we have for all $k$,
\[
  \mathrm{dist}( x_k(t), \; \spt \psi_x ) \ge \frac{1}{3}(L + \sigma_0 t)
\]
where
\[
\sigma_0 := \frac{1}{2} \min_{k=2,\ldots,K} ( c_1^0,\,c_k^0 -
  c_{k-1}^0 ) > 0,
\]
and so
\[
  \norm{u(t)}_{L^{2}(\spt \psi_x)}
  \le C (L+\si_0 t)^{-2} + \norm{\e (t)}_{H^{1/2}(\R)}\ll 1.
\]
The formula $\widehat{\psi_x}(\xi)=e^{-i(x_0+\sigma
t)\xi}\widehat{\zeta'}((b+t)^{\gamma}\xi)$ gives us
\begin{eqnarray*}
  D^{s}\psi_x(x)=\frac{1}{(b+t)^{\gamma(1+s)}}\int 
  e^{i\frac{x-x_0-\sigma t}
  {(b+t)^{\gamma}}\eta}|\eta|^{s}\widehat{\zeta'}(\eta)\,d\eta.
\end{eqnarray*}
Thus
\[
 \||\xi|^{1/2} \widehat{\psi_x} \|_{ L^{1}_\xi} 
  \lesssim (b+t)^{-3\gamma/2}, \quad
  \|\psi\|_{\dot B^{2-2\e}_{\infty,1}}\lesssim  (b+t)^{-2\gamma(1-\e)}.
\]
Similarly,
\[
  \norm{ |\xi|^{1/2} {\cal F}( \sqrt{\psi_x})}_{L^1_\xi} ^2
  \lec  (b+t)^{-2\gamma}.
\]
We can also bound
\[
  \frac \gamma{b+t} |y_k| \le 
  \frac {\si_k} {4(b+t)^\gamma} + 
  \frac {C\gamma^2 y_k^2}{\si_k(b+t)^{2-\gamma}},
\]
and
\[
   \int u^2 \zeta'(y_k) \frac {\gamma^2 y_k^2}
  {\si_k(b+t)^{2-\gamma}}
  \le C (b+t)^{-2+\gamma} [ \norm{\e }_{L^2}^2 + (L+\si_0 t)^{-2}].
\]
Summing the estimates, we get
\begin{align*}
  \frac d{dt} \I_k(t)
  &\le - \frac 12 \int \psi_x |D^{1/2}u|^2
  - \frac {\si_k}4 \int \psi_x u^2
  + C (b+t)^{-3\gamma/2} \norm{u}_{H^{1/2}}^2 \\
  & \quad +  C (b+t)^{-2+\gamma} 
  [ \norm{\e }_{L^2}^2 + (L+\si_0 t)^{-2}].
\end{align*}
Integrating in time and noting $2/3< \gamma < 1$, we get the lemma.
\qed

\subsection{Energy decomposition}
\noindent
{\it Proof of Lemma~\ref{energy}.}
Note $c_k(0) = d_1 + \cdots + d_k$ and
\[
  \sum_{k=1}^K d_k \psi_k
  = \sum_{k=1}^K d_k [\ph_k + \cdots + \ph_K]
  = \sum_{k=1}^K c_k(0) \ph_k.
\]
So
\[
  \G(t) = E(u(t)) + \sum_{k=1}^K d_k \I_k(t) = E(u(t)) + \int_{\mathbb{R}}
  \sum_{k=1}^K \frac 12 c_k(0) \ph_k u^2 \, dx.
\]
Using the decomposition $u = R + \e$ and $R = \sum_{k=1}^K R_k$, we
can decompose $\G(t)$ according to orders in $\e$:
\[
  \G(t) = G_0 + G_1 + \frac12(\e(t), H_K \e(t) ) + \frac 12
  (\e(t),\sumk(c_k(0) - c_k(t))\ph_k \e(t)) - \frac 13 \int_{\R} \e(t)^3,
\]
where $G_0$ denotes terms without $\e$,
\[
  G_0 = E(R) + \frac 12 \int_{\mathbb{R}} \sum_{k=1}^K c_k(0) \ph_k  R^2,
\]
$G_1$ denotes terms linear in $\e$,
\[
  G_1 = \int_{\mathbb{R}} \e[ DR - R^2 + \sum_{k=1}^K c_k(0) \ph_k  R],
\]
and $H_K$ denotes the linear operator
\[
  H_K =  D  - 2R + \sum_{k=1}^K  c_k(t) \ph_k.
\]
We can further decompose
\[
G_0 = \sum_{k=1}^K E(R_k) + \int \sum_{j<k} R_j DR_k - \frac 13
(R^3 - \sum_{k=1}^K R_k^3) + \frac 12 \sum_{k=1}^K c_k(0)  R_k^2
+ \frac 12 \sum_{k=1}^K c_k(0) (\ph_k R^2-R_k^2 ).
\]
Using $DR_k - R_k^2 + c_k(t)R_k = 0$, we have
\[
G_1 = \int \e \bket{
[(\sum_{k=1}^K R_k^2)-R^2]
+ \sum_{k=1}^K \bkt{c_k(0) R_k (\ph_k -1) + c_k(0) \ph_k(R-R_k)
+ (c_k(0) - c_k(t))R_k}}
\]
Note
\[
  \| R^m-\sum_{k=1}^K R_k^m \|_{L^1 \cap L^\infty(\R)}
  \le C L^{-2}, \quad (m=2,3).
\]
\[
  \|R_k (\ph_k -1)\|_{L^2 \cap L^\infty(\R)}+\|\ph_k(R-R_k)\|_{L^2
  \cap L^\infty(\R)} \le C L^{-2}.
\]
Thus
\begin{equation}
  |G_1(t)| \le CL^{-2} + C\sumk |c_k(0) - c_k(t)| 
  \norm{\e}_{L^2}.
\end{equation}
Also, since $R_jDR_k = R_j [c_k(t)R_k - R_k^2]$,
\[
  \|R_jDR_k \|_{L^1 \cap L^\infty(\R)}
  \le C L^{-2}, \quad (j \not = k).
\]
We have
\begin{equation}
  |G_0(t)- \sumk E(R_k) -\sumk c_k(0) N(R_k)| \le CL^{-2}.
\end{equation}
Finally,
\begin{equation}
  \left| \frac 12 (\e,\sumk(c_k(0) - c_k(t))\ph_k \e) 
  - \frac 13 \int \e^3 \right|  \le C \sumk |c_k(0) - c_k(t)| 
  \norm{\e}_{L^2}^2 +  C \norm{\e}_{H^{1/2}}^3,
\end{equation}
completing the proof of Lemma \ref{energy}.
\qed
\medskip

\noindent
{\it Proof of Lemma~\ref{L.Fdiff}.} 
First proof. By energy decomposition around
$Q_{c^0}$, we have for real-valued $\eta$ small in $H^{1/2}$
that
\[
  F_{c^0}(Q_{c^0}+\eta) = F_{c^0}(Q_{c^0}) 
  + \frac 12 (\eta, H^{c^0}\eta) + O(\norm{\eta}_{H^{1/2}}^3).
\]
In particular for $\eta = Q_c - Q_{c^0} $ we get the lemma. In fact
$\eta \sim (c-c^0)\eta_0$ with $\eta_0=\pd_c |_{c=c^0} Q_c$ and $\frac
12 (\eta_0 , H^{c^0}\eta_0)= \frac 12 (\eta_0 , -Q_{c^0})= -\frac 14
\pd_c \int Q_c^2=-\pi/2$.

Second proof. By the scaling property and \eqref{eq:QDQ},
\[
F_{c^0}(Q_c)= c^2 E(Q) + cc^0 N(Q) = -\frac \pi 2 c^2 + cc^0 \pi, \quad
 F_{c^0}(Q_{c^0})  = (c^0)^2 \frac \pi 2.
\]
Thus $F_{c^0}(Q_{c^0}) - F_{c^0}(Q_c) = \frac \pi 2 (c-c^0)^2$.
\qed
\subsection{Quadratic control of $c_k(t)-c_k(0)$}
{\it Proof of Lemma \ref{lem:c}.}
As in the energy expansion above, with 
$u = R + \e$, $R = \sum_j R_j$, and using $( \e, R_j ) \equiv 0$,
and $|x_j(t) - x_k(t)| \geq L/2$ for $j \not= k$,
we have
\[
  |E(u) - \sum_j E(R_j)| \lec 
  L^{-2} + \norm{\e}_{H^{1/2}}^2 + \norm{\e}_{H^{1/2}}^3 .
\]
Now using the conservation of energy, the fact 
$E(R_j) = a c_j^2$, and
$\norm{\e(t)}_{H^{1/2}}\le 1$, we get
\begin{equation}
\label{eq:c1}
  \left|  \sumk \left[ (c_k(t))^2 - (c_k(0))^2 \right] \right|
  \lec L^{-2} 
  + \norm{\e(t)}_{H^{1/2}}^2
  + \norm{\e(0)}_{H^{1/2}}^2.
\end{equation}
Since $\phi_k = \psi_k - \psi_{k+1}$, we have
\[
  {\cal I}_j(t) = \frac{1}{2} \int u^2 \psi_j dx
  = \frac{1}{2} \int u^2 \sum_{k=j}^K \phi_k\,dx
  = \sum_{k=j}^K \int_{\mathbb{R}}\frac{1}{2} \phi_k u^2 dx.
\]
Again using $( \e, R_j ) \equiv 0$,
and $|x_j(t) - x_k(t)| \geq L/2$ for $j \not= k$, 
we see easily that
\[
  \left| \frac{1}{2} \int \phi_k u^2 dx  - N(R_k) \right|
  \lec L^{-2} + \| \e \|_{H^{1/2}}^2.
\] 
So using $N(R_k) =  c_kN(Q_1)=c_k\pi$
and the local monotonicity Lemma~\ref{L2.4}, we get
\begin{equation}\label{eq:c2}
  \de_k(t):=   \sum_{j=k}^K [c_j(t)- c_j(0)] \lec g(t), \quad k=1,\ldots,K,
\end{equation}
where
\[
g(t)= L^{\frac{1}{\gamma}-\frac{3}{2}} + L^{1-\frac{1}{\gamma}}
\sup_{0 \leq \tau \leq t} \norm{\e(\tau)}_{L^{2}}^2 + \| \e(t)
\|_{H^{1/2}}^2 + \| \e(0) \|_{H^{1/2}}^2 .
\]

Denote $\de_{K+1}=0$ and $c_0(0)=0$. Using $|\de| \le -\de + 2\de_+$ for any
$\de \in \R$ and \eqref{eq:c2}, we get
\begin{equation}\label{eq:c3}
\sum_{k=1}^K |\de_k(t)| \lec \sum_{k=1}^K [c_k(0)-c_{k-1}(0)]|\de_k(t) | \le
\sum_{k=1}^K [c_k(0)-c_{k-1}(0)][-\de_k(t) + Cg ].
\end{equation}
By Abel resummation,
\begin{align*}
  -\sum_{k=1}^K [c_k(0)-c_{k-1}(0)]\de_k(t)&=-\sum_{k=1}^K c_k(0)[\de_k(t) -
\de_{k+1}(t)] = \sum_{k=1}^K c_k(0)[c_k(0)-c_k(t)] 
\\
 &= \frac 12 \sumk \bkt{ (c_k(0))^2 -
(c_k(t))^2 } + \frac 12 \sumk  |c_k(t) - c_k(0)|^2.
\end{align*}

Using \eqref{eq:c3}, the above equality and \eqref{eq:c1},  we arrive at
\[
\sumk |\de_k(t)| \lec  g(t)  + \sumk  |c_k(t) - c_k(0)|^2.
\]
Since $|c_k(t) - c_k(0)|\le |\de_k(t)| + |\de_{k+1}(t)|$, we have
\[
\sumk |c_k(t) - c_k(0)|\lec  g(t)  + \sumk  |c_k(t) - c_k(0)|^2.
\] 
By the continuity of $c_k(t)$ and the smallness of $g(t)$, we get Lemma
\ref{lem:c}.  \qed

\subsection{Lower bound for the quadratic form}

We first recall the one-soliton case. Suppose
a function $u(x)$ is a perturbation of $Q_c(x-a)$ of the form
\[
u(x) = Q_c(x-a) + \e (x),
\]
where $\e(x) $ is small in some sense.  Then
\[
(E+cN)(u) = (E+cN)(Q_c) + \frac 12(\e ,H^{c,a} \e ) -
\frac 13 \int \e ^3.
\]
Here $H^{c,a} = D + c - 2Q_c(x-a)$.

\begin{lemma}[\cite{BBSSB}] \label{TH3-1}
Let $H = D+1-2Q$ with $Q(x)= \frac {2}{1+x^2}$.
Its continuous spectrum is $[1,\infty)$. Its eigenvalues are $0,1$, and $\la_\pm =
\frac 12(- 1 \pm\sqrt 5)$, with corresponding normalized
eigenfunctions
\[
\phi_0 = \frac {-4}{\sqrt \pi}\frac {x} {(1+x^2)^2} = \frac {1}{\sqrt
\pi} Q_x,
\]
\[
\phi_1 = \frac {2}{\sqrt \pi} \frac {x(x^2-1)}{(1+x^2)^2} = \frac
{1}{\sqrt \pi}(xQ+Q_x),
\]
\[
\phi_\pm = N_{\pm}\bke{ \frac {1 \pm \sqrt 5}{1+x^2}- \frac 4{(1+x^2)^2}}
= N_{\pm}\bke{(1 + \la_\pm)Q - Q^2}.
\]
Here $N_{\pm} = \frac {1}{\sqrt \pi}(1\pm \frac 2{\sqrt 5})^{1/2}$.
Moreover, there is $\ga_0 \in (0,1)$ so that, if $\e \in H^{1/2}$
satisfies $(\e,Q)=(\e,Q_x)=0$, then
\begin{equation}\label{gap}
\ga_0\norm{\e }_{H^{1/2}}^2 \le (\e , H \e ).  
\end{equation}
\end{lemma}

This lemma, except \eqref{gap}, is due to \cite{BBSSB}. We have
reformulated it in a form convenient to us. To prove Eqn.~\eqref{gap},
decompose $\e = a\phi_- + h$ with $h \perp \phi_-, Q_x$. Thus
\[
(\e,H\e) = \la_- a^2 + (h,Hh) \ge \la_- a^2 + \la_+(h,h) 
= \la_+ (\e,\e) - (\la_+ - \la_-)a^2.
\]
Now decompose $ \phi_{-}= bQ + k$ with $k \perp Q$ and hence
\[
a^2=(\e,\phi_-)^2 = (\e,k)^2 \le (\e,\e)(k,k).
\]
Thus
\[
(\e,H\e) \ge \ga (\e,\e), \quad \ga= \la_+ - (\la_+-\la_-)(k,k)
\]
One can compute $(k,k)=\frac 12 - \frac 1{\sqrt 5}$ and $\ga =\frac 12$,
and Eqn.~\eqref{gap} follows with $\ga_0=\frac 19$.

We can rescale \eqref{gap} and get the following: Let $R(x)=Q_c(x-a)$.
If $\e \in H^{1/2}(\R)$ satisfies $(\e,R)=(\e,R_x)=0$, then
\begin{equation}
\ga_0(\e,(D+c)\e)  \le (\e , (D+c-2R) \e ).  
\end{equation}

{\it Proof of Lemma~\ref{TH3-2}.} This is a time-independent statement
and everything is evaluated at $t$, e.g., $c_k=c_k(t)$.
Let $\chi(x)$ be a nonnegative smooth function supported in 
$|x|\le 2$, $\chi(x)=1$ for $|x|\le 1$, and $\chi^2(x) \le 1/2$ 
if and only if $|x| \ge 3/2$. 
Let $\chi_k(x) = \chi(\frac{x-x_k}{L_2/16})$.  In
particular $\ph_k(x)=1$ when $\chi_k(x) \not = 0$, and $\ph_k(x) \ge
2\chi_k^2(x)$ when $\chi_k^2(x) \le 1/2$.
Decompose
\begin{align*}
(\e, H_K\e) =& \sumk  ( \chi_k \e, 
(D+c_k - 2R_k) (\chi_k \e))   
\\
&+(\e,D\e) - \sumk (\chi_k \e, D (\chi_k \e))
\\
&+ \sumk c_k (\e,  (\ph_k - \chi_k^2)\e)
\\
&+ (\e, -\sumk 2R_k(1- \chi_k^2)\e)
\\
=:&\ I_1+I_2+I_3+I_4.
\end{align*}
It follows from Lemma \ref{TH3-1} that
\begin{equation}\label{eq3-1}
I_1 \ge \sumk \ga_0 (\chi_k \e, (D+c_1) (\chi_k \e)) .
\end{equation}
By Lemma \ref{L2.2}
\[
\left| (\chi_k \e, D (\chi_k \e))^{1/2}- \norm{\chi_k D^{1/2} \e
 }_{L^2}\right | \le \norm{[D^{1/2},\chi_k ] \e}_{L^2} \le 
\norm{|\xi|^{1/2} \hat \chi_k(\xi)}_{L^1(d\xi)}\norm{\e}_{L^2}.
\]
By definition of $\chi_k$,
\[
\norm{|\xi|^{1/2} \hat \chi_k(\xi)}_{L^1(d\xi)} \le C L_2^{-1/2}.
\]
Thus
\[
  I_1 \ge \text{RHS of }\eqref{eq3-1}
  \ge \frac {\ga_0}2  \sumk \norm{\chi_k D^{1/2} \e} _{L^2}^2 -
  CL_2^{-1}  \norm{\e} _{L^2}^2
  + \ga_0 c_1 \sumk\norm{\chi_k \e}_{L^2}^2 ,
\]
and
\begin{align*}
I_2 &\ge (\e, D\e) - (1 + \frac {\ga_0}4) \sumk \norm{\chi_k D^{1/2} \e}
_{L^2}^2 - CL_2^{-1} \norm{\e} _{L^2}^2
\\
&= \int [1-(1+\frac{\ga_0}4)\sumk\chi_k^2] |D^{1/2} \e|^2  
- CL_2^{-1} \norm{\e} _{L^2}^2.
\end{align*}
We also have 
\[
  I_3 \ge \sumk c_1 (\e, \ph_k 1(\chi_k^2 \le 1/2)\e) 
  = c_1 \int _{ \sum_k \chi_k ^2\le 1/2} \e^2,
\]
\[
  |I_4| \le C L_2^{-2} (\e,\e).
\]
Summing up, we have
\[
  (\e, H_K\e) \ge \frac {\ga_0}4 (\e,D\e) - CL_2^{-1}\|\e\|_2^2
  + \ga_0 c_1 \sumk\norm{\chi_k \e}_{L^2}^2 + 
  c_1 \int _{ \sum_k \chi_k ^2\le 1/2} \e^2
\]
which is greater than $\frac {\ga_0} 4(\e,(D+c_1)\e)$ if $L_2$ is
sufficiently large.
\qed

\section*{Acknowledgments}
The research of Gustafson and Tsai is partly
supported by NSERC grants (Canada).
The research of Takaoka is partly
supported by a JSPS grant (Japan).


\end{document}